\documentclass[a4paper, 12pt]{article}

\usepackage[english]{babel}
\usepackage{authblk}

\usepackage{amsfonts,amsmath,amsthm,mathrsfs,amscd,amssymb}
\usepackage{cite}

\newtheorem{theorem}{Theorem}

\newtheorem{prop}{Proposition}

\providecommand{\keywords}[1]
{
  \small	
  \textit{Key Words:~} #1
}

\providecommand{\MSC}[1]
{
  \small	
  \textit{Mathematics Subject Classification 2020:~} #1
}

\newcommand*{\cm}[1]{\mathscr{#1}}

\begin{document}

\title{On Khinchin's theorem about the special role of the Gaussian distribution \footnote{This is a preprint of the work accepted for publication in the Journal of Contemporary Mathematical Analysis \textbf{58} (6), 2023, the copyright holder indicated in the Journal}}

\author[1]{Linda A. Khachatryan}

\affil[1]{\small Institute of Mathematics, NAS RA, {\it linda@instmath.sci.am}}

\date{}

\maketitle

\begin{abstract}
The purpose of this note is to recall one remarkable theorem of Khinchin about the special role of the Gaussian distribution. This theorem allows us to give a new interpretation of the Lindeberg condition: it guarantees the uniform integrability of the squares of normed sums of random variables and, thus, the passage to the limit under the expectation sign. The latter provides a simple proof of the central limit theorem for independent random variables.
\end{abstract}


\MSC{60F05, 60G50}

\bigskip

\noindent \textbf{1.} Let $\{\xi_{n,j}\} = \{\xi_{n,j}, 1 \le j \le k_n, n \ge 1\}$, $k_n \to \infty$ as $n \to \infty$, be a triangular array (double sequence) of independent in each row random variables on a probability space $(X, \cm{B}, P)$. For the sake of simplicity, we always assume that $E\xi_{n,j} = 0$ for all $j$ and $n$. For any $n \ge 1$, denote $S_n = \sum \limits_{j=1}^{k_n} \xi_{n,j}$, and let $DS_n$ be its variance. The Gaussian (normal) distribution function with parameters $a$ and $\sigma^2$, $a,\sigma \in \mathbb{R}$, $\sigma > 0$, is defined by
$$
\Phi_{a,\sigma^2} (x) = \frac{1}{\sigma \sqrt{2\pi}} \int \limits_{-\infty}^x \exp\left\{ - \frac{(t-a)^2}{2 \sigma^2} \right\} dt, \qquad x \in \mathbb{R}.
$$

Khinchin~\cite{Khin} (translation into English can be found in~\cite{RogMain}) noted that the Gauss law, as a limiting law for sums of independent random variables, has a very special role that distinguishes it from all infinitely divisible laws. Namely, we arrive at the Gauss law in all cases when the limiting negligibility of the components of the sum of terms under study reaches a sufficiently strong degree; and this happens completely independently of the special properties of the laws of distribution of these terms.

The condition of asymptotic infinitesimality (or, equivalently, limiting negligibility) on the summands $\xi_{n,j}$, in the general case, is formulated as the condition that for any $\varepsilon >0$, probability of the inequality $\vert \xi_{n,j} \vert \ge \varepsilon$ tends to zero uniformly in $j$ as $n \to \infty$:
\begin{equation}
\label{infsmall1}
\max \limits_{1 \le j \le k_n} P \left(\vert \xi_{n,j} \vert \ge \varepsilon \right) \to 0 \quad \text{as } n \to \infty.
\end{equation}
Khinchin showed (Theorem 42 in~\cite{Khin}) that if we assume that not only this probability but the probability that all $\vert \xi_{n,j} \vert$, $1 \le j \le k_n$, are greater than $\varepsilon$ tends to zero as $n \to \infty$, that is,
\begin{equation}
\label{infsmall2}
P \left( \max \limits_{1 \le j \le k_n} \vert \xi_{n,j} \vert \ge \varepsilon \right) \to 0 \quad \text{as } n \to \infty,
\end{equation}
then the only possible limiting law for normed row sums is the Gauss law.

\begin{theorem}[Khinchin]
Let $\{\xi_{n,j}\}$ be a double sequence of independent in each row random variables. If a limiting non-degenerate distribution for the sums $S_n$ exists, then for it to be Gaussian, it is necessary and sufficient that for any $\varepsilon>0$, random variables $\{\xi_{n,j}\}$ satisfy~\eqref{infsmall2}.
\end{theorem}

Since condition~\eqref{infsmall2} represents only a somewhat strengthened requirement~\eqref{infsmall1} for the limiting negligibility of summands and does not contain any special assumptions about the nature of the laws of distribution of summands, the above result characterizes the Gauss law as, in a certain sense, a universal limiting law for sums of independent random variables and justifies the exclusive place given to this law in classical studies.

In the bibliographical notes~\cite{Khin}, Khinchin mention, that a more general result was obtained by L\'{e}vy in~\cite{Levy}, however, Khinchin was not able to find the proof of the latter based on the sketch suggested by L\'{e}vy. Khinchin's proof is based on the direct investigation of the characteristic functions of summands. Another (shorter) proof, based on the L\'{e}vy--Khinchin formula for the decomposition of characteristic functions, was suggested by Gnedenko~\cite{Gned}. The latter can be found in the book~\cite{GnedKolm} by Gnedenko and Kolmogorov (see Theorem 1 on p. 126).

Under conditions of the Khinchin theorem, the limiting distribution (in the case of centered random summands) is $\Phi_{0,\sigma^2}$ with some parameter $\sigma^2$. Since the limiting law is non-degenerate, $\sigma^2>0$. Note that Khinchin do not impose any restriction on the second moments of the summands, which is a minimal condition on the moments in the central limit theorem (CLT). We say that for a sequence $\{\xi_{n,j}\}$ of (centered) random variables, the CLT holds if
$$
\lim \limits_{n \to \infty} P\left( \frac{S_n}{\sqrt{DS_n}} \le x \right) = \Phi_{0,1}(x), \qquad x \in \mathbb{R}.
$$

The Khinchin theorem cannot be considered a CLT since the limiting Gaussian distribution is not necessarily the standard one with $\sigma^2 = 1$. However, if we impose the uniform integrability condition on the squares of normed row sums, we will be able to prove CLT based on the Khinchin result.

We will need the following statements  (see, for example, Lemma 1 on p. 322 and Theorem 5 on p. 189 in~\cite{Shir}). Under convergence of distribution functions we understand convergence in general, i.e., at each point of continuity of the limiting distribution function.

\begin{prop}
\label{proposition}
Let $\{F_n\} = \{F_n, n \ge 1\}$ be a sequence of distribution functions. Suppose that any convergent subsequence $\{F_{n'}\}$ of $\{F_n\}$, $\{n'\} \subset \{n\}$, converges to the same distribution function $F$. Then the sequence $\{F_n\}$ converges to $F$ as well.
\end{prop}
\begin{proof}
Let $\cm{X}_F$ be the set of continuity points of the distribution function $F$. Fix some $x \in \cm{X}_F$ and assume that $F_n(x)$ does not converge to $F(x)$. Then there exists $\varepsilon > 0$ and an infinite sequence $\{n'\}$ of natural numbers such that
\begin{equation}
\label{Fn->F}
\left| F_{n'}(x) - F(x) \right| > \varepsilon.
\end{equation}
By the Helly theorem, from the sequence $\{F_{n'}\}$, one can select a convergent subsequence $\{F_{n''}\}$, and let generalized distribution function $G$ be its limit. By the hypothesis of the proposition, $G = F$, and thus, $F_{n''}(x) \to F(x)$ as $n \to \infty$, which contradicts with~\eqref{Fn->F}. This completes the proof.
\end{proof}

We remind that a family of random variables $\{\eta_n, n \ge 1\}$ is uniformly integrable if
$$
\sup \limits_{n} \int \limits_{\vert \eta_n \vert > C} \vert \eta_n \vert dP \to 0 \qquad \text{as } C \to \infty.
$$

\begin{theorem}
\label{th-unif-integr}
Let $\eta_n$, $n \ge 1$ be a sequence of positive random variables with $E\eta_n < \infty$ such that $\eta_n \to \eta$ as $n \to \infty$. Then $E\eta_n \to E\eta < \infty$ as $n \to \infty$ if and only if the family $\{\eta_n, n \ge 1\}$ is uniformly integrable.
\end{theorem}

Now we present the following version of the CLT for independent random variables.

\begin{theorem}
\label{th-main}
Let $\{\xi_{n,j}\}$ be a double sequence of independent in each row random variables such that $E\xi_{n,j}^2 < \infty$, $1 \le j \le k_n$, $n \ge 1$. If random variables $\{\xi_{n,j}/\sqrt{DS_n}, 1 \le j \le k_n, n \ge 1\}$ satisfy condition~\eqref{infsmall2} and the squares of normed row sums $\{S_n^2/DS_n, n \ge 1\}$ are uniformly integrable, then for the sequence $\{\xi_{n,j}\}$, the CLT holds.
\end{theorem}
\begin{proof}
Let $F_n$ be the distribution function of $S_n/ \sqrt{DS_n}$, $n \ge 1$. Then
\begin{equation}
\label{D=1}
\int \limits_{\mathbb{R}} x^2 dF_n(x) = E\left( \frac{S_n^2}{DS_n} \right) = 1, \qquad n \ge 1.
\end{equation}
Further, let $\{F_{n'}\}$, $\{n'\} \subset \{n\}$ be some convergent subsequence of the sequence $\{F_n\}$. Due to the Khincin theorem, $F_{n'}(x) \to \Phi_{0,\sigma^2}(x)$ as $n' \to \infty$ for any $x \in \mathbb{R}$ and some $\sigma >0$ if random variables $\{\xi_{n,j}/\sqrt{DS_n}\}$ satisfy condition~\eqref{infsmall2}. Since $\{S_n^2/DS_n, n \ge 1\}$ are uniformly integrable, due to Theorem~\ref{th-unif-integr}, we can pass to the limit under the expectation sign, and thus,
$$
\lim \limits_{n' \to \infty} \int \limits_{\mathbb{R}} x^2 dF_{n'}(x) = \int \limits_{\mathbb{R}} x^2 \lim \limits_{n' \to \infty} dF_{n'}(x) = \int \limits_{\mathbb{R}} x^2 d\Phi_{0,\sigma^2}(x).
$$
Tacking into account~\eqref{D=1}, we conclude
$$
\int \limits_{\mathbb{R}} x^2 d\Phi_{0,\sigma^2}(x) = 1,
$$
that is, the parameter $\sigma^2$ in the limiting Gaussian distribution is equal to one. Thus, from the uniform integrability of $\{S_n^2/DS_n, n \ge 1\}$, it follows that $F_{n'}(x) \to \Phi_{0,1}(x)$ as $n' \to \infty$ for any $x \in \mathbb{R}$.

Thereby, any convergent subsequence $\{F_{n'}\}$ of the distribution functions of the normed row sums $S_n/ \sqrt{DS_n}$ converge to the same limiting distribution $\Phi_{0,1}$. Hence, by Proposition~\ref{proposition}, the sequence $\{F_{n}\}$ converges to $\Phi_{0,1}$ as well. Therefore, for random variables $\{\xi_{n,j}\}$, the CLT holds.
\end{proof}

%Let us note, that if random variables $\{\xi_{n,j}/\sqrt{DS_n}, 1 \le j \le k_n, n \ge 1\}$ satisfy condition~\eqref{infsmall2}, then from the validity of the CLT for the sequence $\{\xi_{n,j}\}$ it follows the uniform integrability of the squares of normed row sums $\{S_n^2/DS_n, n \ge 1\}$.

\noindent \textbf{2.} Theorem~\ref{th-main} allows us to give the new probabilistic interpretation of the Lindeberg condition. We will show that from the Lindeberg condition, the uniform integrability of the squares of normed sums of random variables follows, which, in its turn, allows passage to the limit under the expectation sign, and thus, guaranties $\sigma^2 = 1$ in the limiting Gaussian distribution in the Khinchin theorem.

The double sequence $\{\xi_{n,j}\}$ of random variables satisfies the Lindeberg condition if for any $\varepsilon > 0$,
\begin{equation}
\label{Lind}
\frac{1}{DS_n} \sum \limits_{j=1}^{k_n} \int \limits_{\{\vert \xi_{n,j} \vert > \varepsilon \sqrt{DS_n}\}} \xi_{n,j}^2 dP \to 0 \qquad \text{as } n \to \infty.
\end{equation}
The classical interpretation of the Lindeberg condition is that if a sequence $\{\xi_{n,j}\}$ of random variables satisfies~\eqref{Lind}, then its elements are asymptotically infinitesimal uniformly in each row, that is, relation~\eqref{infsmall1} holds. Billingsley (see p. 90 in~\cite{Bill}) noted that from the Lindeberg condition, the uniform integrability of the squares of normed sums follows as well.

\begin{prop}
\label{prop-Lind-interpr}
Let $\{\xi_{n,j}\}$ be a double sequence of (centered) independent random variables with finite second moments. If $\{\xi_{n,j}\}$ satisfies the Lindeberg condition~\eqref{Lind}, then squares of normed row sums $\{S_n / \sqrt{DS_n}, n \ge 1\}$, are uniformly integrable.
\end{prop}
\begin{proof}
The statement follows from inequality (12.20) in~\cite{Bill}, according to which for any $n \ge 1$ and $C > 0$, one has
$$
\int \limits_{\{S_n^2 \ge C DS_n\}} \frac{S_n^2}{DS_n} dP \le K \left( \frac{1}{C} + \frac{1}{DS_n} \sum \limits_{j=1}^{k_n} \int \limits_{\{\vert \xi_{n,j} \vert \ge \frac{1}{4} C DS_n\}} \xi_{n,j}^2 dP \right)
$$
where $K$ is some universal constant. By~\eqref{Lind}, for any $C>0$, there exists $n_0 = n_0(C) > 1$ such that
$$
\frac{1}{DS_n} \sum \limits_{j=1}^{k_n} \int \limits_{\{\vert \xi_{n,j} \vert \ge \frac{1}{4} C DS_n\}} \le \frac{1}{C} \qquad \text{for any } n > n_0.
$$
Thus,
$$
\sup \limits_n \int \limits_{\{S_n^2 \ge C DS_n\}} \frac{S_n^2}{DS_n} dP \le K \left( \frac{2}{C} + \sup \limits_{1 \le m \le n_0(C)} \frac{1}{DS_m} \sum \limits_{j=1}^{k_m} \int \limits_{\{\vert \xi_{m,j} \vert \ge \frac{1}{4} C DS_m\}} \xi_{m,j}^2 dP \right),
$$
and hence,
$$
\sup \limits_n \int \limits_{\{S_n^2 \ge C DS_n\}} \frac{S_n^2}{DS_n} dP \to 0 \qquad \text{as } C \to \infty.
$$
\end{proof}

The statement above reveals the true essence of the Lindeberg condition. Since uniform integrability condition is the necessary and sufficient condition for taking limit under the expectation sign, we conclude that the Lindeberg condition is one of conditions under which the limiting Gaussian distribution in the Khinchin theorem is the standard one. Tacking into account this fact, we provide the new proof of the well-known L\'{e}vy-Lindeberg theorem.

\begin{theorem}[L\'{e}vy-Lindeberg]
\label{th-L-L}
Let $\{\xi_{n,j}\}$ be a double sequence of independent in each row random variables such that $E\xi_{n,j} = 0$, $0< E\xi_{n,j}^2 < \infty$, $1 \le j \le k_n$, $n \ge 1$. If random variables $\{\xi_{n,j}\}$ satisfy Lindeberg condition~\eqref{Lind}, then the CLT holds.
\end{theorem}
\begin{proof}
First note that random variables $\{\xi_{n,j}/\sqrt{DS_n}, 1 \le j \le k_n, n \ge 1\}$ satisfy condition~\eqref{infsmall2}, since for any $\varepsilon>0$, we have
$$
P \left(\max \limits_{1 \le j \le k_n}  \vert \xi_{n,j} \vert \ge \varepsilon \sqrt{DS_n}\right) \le \sum \limits_{j=1}^{k_n} P(\vert \xi_{n,j} \vert \ge \varepsilon \sqrt{DS_n}) \le \frac{1}{\varepsilon^2 DS_n} \sum \limits_{j=1}^{k_n}  \int \limits_{\{\vert \xi_{n,j} \vert > \varepsilon \sqrt{DS_n}\}} \xi_{n,j}^2 dP \to 0
$$
as $n \to \infty$. Further, by Proposition~\ref{prop-Lind-interpr}, random variables $\{S_n^2/DS_n, n \ge 1\}$, are uniformly integrable. Thus, by Theorem~\ref{th-main}, for $\{\xi_{n,j}\}$ the CLT holds.
\end{proof}

\noindent \textbf{3.} Let us illustrate the application of Theorem~\ref{th-main} in the case of independent identically distributed (i.i.d.) random variables. Namely, we will use this theorem to prove the following classical result.

\begin{theorem}[L\`{e}vy--Khinchin]
	Let $\{\eta_n, n \ge 1\}$ be a sequence of i.i.d. random variables such that $E\eta_1=0$ and $D\eta_1 = \sigma_0^2 \le \infty$. Then for $\{\eta_n\}$, the CLT holds.
\end{theorem}
\begin{proof}
	Consider the double array $\{\xi_{n,j}, 1 \le j \le n, n \ge 1\}$ of random variables $\xi_{n,j} = \dfrac{\eta_j}{\sigma_0 \sqrt n}$. It is not difficult to check, that random variables $\{\xi_{n,j}\}$ satisfy condition~\eqref{infsmall2}. Further, put $\hat S_n = \sum \limits_{j=1}^n \xi_{n,j} = \dfrac{1}{\sigma_0 \sqrt n} \sum \limits_{j=1}^n \eta_j$, then $D \hat S_n = 1$. With application of inequality (12.19) in~\cite{Bill}, for any $C>0$, we can write
	\begin{eqnarray*}
		P(\hat S_n^2 \ge C^2) & \le & \max \limits_{1 \le j \le n} P \left(\hat S_j^2 \ge C^2 \right) \le P\left( \max \limits_{1 \le j \le n} \vert \hat S_j \vert \ge C \right) \le \\
		\\
		& \le & K \left( \dfrac{1}{C^4} + \displaystyle\dfrac{1}{C^2} \sum \limits_{k=1}^n \int \limits_{\{\vert \xi_{n,j} \vert > \frac{1}{4} C\}} \xi_{n,j}^2 dP \right),
	\end{eqnarray*}
	where $K$ is some positive constant. Further, applying equality (3) on p. 223 in~\cite{Bill}, we can write
	\begin{eqnarray*}
		\displaystyle\int \limits_{\hat S_n^2 \ge C} \hat S_n^2 dP & = & C P \left(\hat S_n^2 \ge C \right) + \int \limits_{C}^\infty P\left( \hat S_n^2 \ge t\right) dt \le \\
		\\
		& \le & C \left( \dfrac{K}{C^4}  + \dfrac{K}{C^2\sigma_0^2} \displaystyle\int \limits_{\{\vert \eta_1 \vert > C \sigma_0/4\}} \eta_1^2 dP\right) + \displaystyle\int \limits_{C}^\infty \left( \dfrac{K}{t^4}  + \dfrac{K}{t^2\sigma_0^2} \displaystyle\int \limits_{\{\vert \eta_1 \vert > t \sigma_0/4\}} \eta_1^2 dP\right) dt = \\
		\\
		& = & K \left( \dfrac{1}{C^3} + \displaystyle\int \limits_{C}^\infty \dfrac{dt}{t^4} + \dfrac{1}{\sigma_0^2} \left( \dfrac{1}{C} + \int \limits_{C}^\infty \dfrac{dt}{t^2} \right) \int \limits_{\{\vert \eta_1 \vert > t \sigma_0/4\}} \eta_1^2 dP   \right).
	\end{eqnarray*}
From here it follows that random variables $\{\hat S_n, n \ge 1\}$ are uniformly integrable. Hence, by Theorem~\ref{th-main}, for random variables $\{\xi_{n,j}\}$ the CLT holds. It remains to note that $P \left(\hat S_n \le x \right) = P\left(\dfrac{S_n}{\sqrt{DS_n}} \le x \right)$, $x \in \mathbb{R}^d$.
\end{proof}

We see, that checking the uniform integrability of the squares of normmed sums directly requires some effort. At the same time, the Lindeberg condition for i.i.d. random variables can be checked quite simply. That is why it is preferable in applications.\\

\noindent \textbf{Acknowledgments.} The author is grateful to Boris S. Nahapetian for his attention to the work and useful discussions. Also, the author wants to express her gratitude to Lutz Mattner for his comments on the preprint of this work, including the reference to the book~\cite{RogMain}, and to an anonymous reviewer for valuable remarks.

\end{document}